\documentclass{article}
\usepackage{amsfonts,amssymb,amsmath,amsgen,amsopn,amsbsy,theorem,graphicx,epsfig}
\usepackage{amscd,bezier,latexsym,mathrsfs,enumerate,multirow}
\usepackage[utf8]{inputenc}\usepackage[english]{babel}
\usepackage[dvipsnames]{xcolor}
\usepackage[pagewise]{lineno}
\newcommand{\bc}{\begin{center}}
\newcommand{\ec}{\end{center}}

\numberwithin{equation}{section}

\newtheorem{theorem}{Theorem}[section]

\newtheorem{definition}[theorem]{Definition}

\title{IN MEMORIAM\\
 CEM TEZER -- (1955-2020) }
\author{Fatma Muazzez Şimşir \\
Selçuk University}
\date{}  

\begin{document}
\maketitle

\begin{abstract}Cem Tezer was a fastidious, meticulous, highly idiosyncratic and versatile scientist. Without him 
Turkish community of mathematics would be incomplete. Our sense of gratitude for his work in various areas of mathematics, history of sciences, literature, music and his encouragement to do mathematics for only its beauty was hardly unique and even unusual. After he passed away on 27 February 2020, while working actively at Middle East Technical University, the number of colleagues and former students described the ways in which their studies and indeed their view towards mathematics had been transformed by having known him might have surprised only those who had never met him. In this article not only, his contributions to mathematics will be classified and summarized but also his unique and distinguished personality as a mathematician will be emphasized.

{\bf Keywords: }{Dynamical systems and ergodic theory, general topology, classical and differential geometry, global analysis, history of sciences and mathematics, literature, music.}
\end{abstract}

\section{Introduction}

\begin{figure}[ht]
\centering
\includegraphics{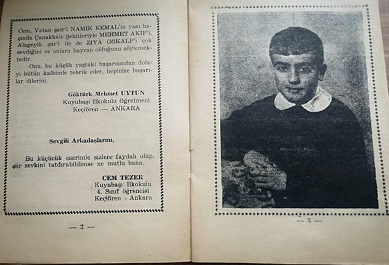}
\caption{Deniz Feneri, Şiirler.}
\label{fig1}
\end{figure}

Cem Tezer passed away in Dikmen, Ankara on 27 February 2020, from a sudden heart attack in the arms of his brother Uğur Tezer.  With the passing away of Cem Tezer  Turkish Mathematical community not only lost a distinguished researcher and colleague but also is left without an eloquent orator. He was born in Bilecik on March 19, 1955 as the first son of Orhan and Sevin\c{c} Tezer. Orhan Tezer was a Civil Engineer with a MSc. degree and Sevinç Tezer left Ankara University Faculty of Languages History and Geography (DTCF) while she was studying philosophy. The family settled down in several different cities of Anatolia due to his father's employment as a Civil Engineer. 
He studied first two years of the primary school at İstiklal İlkokulu (1960-1962) in Merzifon. This little town left him deep with impressions and lifelong friendships.
After Merzifon, they settled down to Keçiören, Ankara (1962-1965) and he finished the primary school at Kuyubşı İlkokulu where his teacher 
G\"okt\"urk Mehmet Uytun discovered his talent for poetry and encouraged him to publish his poems in a booklet entitled as "{\em{Deniz Feneri}}, Şiirler", \cite{1}. For the foreword of this book with a lovely picture of Cem with school uniform, see Figure \ref{fig1}. 
After studying the first year of the secondary school in Kalaba Ortaokulu, they moved to Elazığ (1965-1967) where he took accordion courses, met with his dear friend Erkan Oğur and   finished his  secondary school education. In this period of life he engaged in music with a great love and at the same time he developed a passion 
towards chemistry, even he built a laboratory in his house where he carried out experiments. 
 
In 1965, Tezer family returned to Ankara and Cem Tezer graduated from high school in Ankara at Atatürk Lisesi. 
During those years his interest in mathematics dominated chemistry. By the influence and support  of his father 
he started university at the Civil Engineering Faculty of İstanbul Technical University in 1971. Somehow his father persuaded him that
he could learn Mathematics as an engineering student but after a few months he recognized the truth and decided to leave. 
As an influential character his father  advised him to study abroad if he insisted on doing mathematics. After then, he went to London, started to learn English and
prepared for the exams of Cambridge University. 

During his studies at Cambridge University (1973-1976), his love in music revived again and he intensely dealt with European chamber music. First,
he learned the  block flute, played especially pieces from old English masters and baroque era and then bought a flute. 
I still remember vividly the joy and the pleasure of listening to Cem Tezer rehearsing the flute in the Department of Mathematics at Middle East Technical University occasionally late in the afternoon as a bachelor student. As a young bachelor student in Cambridge he was mainly interested in the foundations of mathematics and philosophy yet
a remarkable amount of courses were related to physics including electromagnetism, thermodynamics and quantum mechanics. After receiving his B.Sc. degree in mathematics at Cambridge University (1976) he started at his doctoral studies in Heidelberg and earned his Ph.D. degree at Ruprecht-Karls-Universität Heidelberg (1984). His thesis is entitled as {\em "Generalised Solenoids of the Dantzig Type”}, \cite{2}. Professor Tezer later published the results of his thesis in the Israel Journal of Mathematics in 1987  under the title "The shift on the inverse limit of a covering  projection", \cite{3}. 

After earning his Ph.D. he returned to Turkey  he decided to work at Middle East Technical University, Ankara where he would spend the next almost four decades where he was involved extensively in research and teaching. He also had a book chapter that is very carefully  written about the foundations of algebraic topology,\cite[Ch.\ 24, p. 779-848]{26}.

In the following one may find his M.Sc. and Ph.D. students with the titles of their thesis in chronological order.

\begin{itemize}
\item M. Işık Güzelgöz, M.Sc., Affine Conformal Vector Fields, (1996).
\item Nilüfer Koldan, M.Sc., Prevalence Of Almost Inner Automorphisms In Isospectral Deformations Of Riemannian Two Step Nilmanifolds,(2001).
\item Fatma Muazzez Şimşir, M.Sc., On Two Instances Of Spectral Rigidity, (2001).
\item Semra Taşkın(Pamuk), M.Sc., Connection Preserving Conformal Diffeomorphisms Of Spheres, (2002).
\item Selma Yıldırım, M.Sc.,  Magnetic Spherical Pendulum, (2003).
\item Fatma Muazzez Şimşir, Ph.D., Conformal Vector Fields With Respect To The Sasaki Metric (2005).
\end{itemize}

\subsection*{Glimpses about him as a polymath}

The influence on Cem Tezer of the rich cultural, historical and academic heritage he observed in these two outstanding institutions, was the main source of inspiration in his academic and intellectual work. This overwhelming influence led him to a principle which he employed at each stage of his mathematical research; "study only the classical works of the true masters of the field". He was keen to avoid "reader friendly" mathematical expositions. He was also critical of fashionable techniques most of which did fade away very rapidly without leading to any substantial progress. Cem Tezer adored mathematics as a human endeavour. He did not try to comply with the dominating academic attitude formulated in the cliché "publish or perish"; on the contrary, his attitude was that of the masters – few but ripe. His research articles span four basic fields, namely topological dynamics, differential geometry, mechanics and elementary geometry. I always had the impression that among these fields, elementary geometry was his favourite; he quite often emphasized the neatness, simplicity and the beauty of the ideas and the techniques in elementary geometry.
The ideas, the methods and the results in each and every research article of Cem Tezer deserve a careful investigation, which I believe will lead to further interesting research, \cite{25}. Those were the late eighties and early nineties that he was very keen on working with Hüseyin Demir on elementary geometry and they published
articles and answered several questions from the question series asked by American Mathematical Monthly and Mathematical Magazine, \cite{16}--\cite{22}.  
Matematik Dünyası was also another opportunity for him to share his original, beautiful proof techniques on Euclidean geometry.\footnote{Özlük, Ö., Şahin, A.  ve Tezer,C. Pisagor Teoreminin Çeşitli İspatları, Matematik Dünyası 1, 3(1991) pp.  6-9. \\
Menelaus ve Ceva Teoremleri ve Düzlem Geometride Açılar ve Ölçüleri\\ 
https://www.matematikdunyasi.org/article/menealus-ve-ceva-teoremleri/.}

\begin{figure}[h!]
\centering
\includegraphics{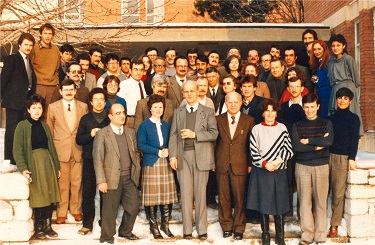}
\caption{Backyard of Department of Mathematics, METU with Hüseyin Demir in the center.}
\label{fig2}
\end{figure}

Cem Tezer had great interest in the history of mathematics pertaining to the late Ottoman era. Combining his deep mathematical insight with his competence in the language and the script in which the articles were written then, Cem Tezer contributed significantly to the studies on the mathematics of that period. On several occasions he was invited to  DTCF (Faculty of Language History and Geography) at Ankara University, one of the leading departments of history in Turkey, to teach graduate courses in history of sciences and mathematics. Başhoca İshak Efendi and Vidinli Hüseyin Tevfik Paşa both had a central role in his research, \cite{11}, \cite{12}. He had also an  article about Hüseyin Demir, \cite{13}.

His speech was urbane and sophisticated, and his command of Turkish, English, German and French  was eloquent in conversation as it was on the page. He was able to read and write in Persian. He was more than capable of expressing sharp and accurate judgements, and independent of the  language he spoke,
his discourse was in terms so carefully chosen and so persuasive that dissent would have seemed superfluous.

Literature and music were two other endeavours which fascinated Cem Tezer. In his poems, mainly in Turkish, a few in English and in unpublished articles on Turkish literature we again observe how keen he was to prefer the classical over the fashionable. In his poems in Turkish, Cem Tezer almost invariably followed the classical tradition, using the rich late Ottoman era language and style. He derived his inspiration from the masterpieces of Divan Edebiyatı. Cem Tezer was about to complete a trilogy on three of the most influential and controversial Turkish poets; Yahya Kemal Beyatlı, Tevfik Fikret and Mehmet Akif Ersoy. I believe that even though the finishing touches of Cem Tezer are lacking, these articles are of great importance as terse comparative essays on various social, cultural and literary attitudes prevailing in the Ottoman Empire and then in the Turkish Republic in early 20th century, \cite{25}. 

Cem Tezer had several poems in different folders and he shared with close friends. It was my honour that he also shared one of those with me called "Bazı Bir Sırça Muamma", \cite{14}. Among his poems, those he quintupled "tahmis etmek"  the poems called "Gazel"'s of masters of poetry  such as N$\hat{a}$bi, Yahya Kemal, Turgut Uyar and Attila İlhan plays a central role of his work reflecting his style.\footnote{ The last quintet of Yahya Kemal'in Veda Gazelini Tahmis
\\Gösterdi ehil ehlini tahmis-i gazelde !\\ 
Eyyam-ı Kemal-veş tükenip sonsuz emelde\\
Gül solsa kalıp s$\hat{ı}$m, siyah kalsa da elde\\
\em tekrar mülaki oluruz bezm-i ezelde:\\
Evvel giden ahbaba selam olsun erenler ! }

\begin{figure}[h!]
\centering
\includegraphics{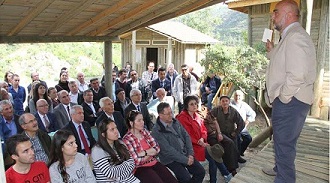}
\caption{Çakılarası Mathematics Village, Opening Talk, 2016.}
\label{fig3}
\end{figure}

The last but not the least, he was a distinguished teacher with remarkable energy and enthusiasm.\footnote{Math 373 - Geometry by Cem Tezer:https://ocw.metu.edu.tr/course/view.php?id=311} As far as I know he was about to finish a textbook  in differential
geometry of curves and surfaces in Turkish. However, this time his meticulousness deprived  us from a rigorously written  
textbook which will be unique in terms of both its language and content.  In written expositions, numerous conference and seminar talks \footnote{Türklerin Matematiğe Katkıları: https://www.youtube.com/playlist?list=PLBaOxR9yqKjVDwLXXqM0oyN9gTNRJ9o5y} and in private discussions he generously shared intriguing ideas, questions in mathematics, history, literature, languages, music. He enjoyed every aspect of civilisation like
concert halls, elegant cafés and restaurants, museums, libraries, book stores yet once the subject is mathematics contrary to his city motivated joys he loved to be in
Çakılarası Mathematics Village\footnote{https://www.cakilarasimatematikkoyu.com/} where he gave several lecture series. These lecture series had a great diversity including "Gamma and Zeta functions","Theory of special relativity", "Euclidean geometry",  "Non-Euclidean geometry", "Characterisation of closed surfaces",  see Figure \ref{fig3}.

When he deceased on February 2020, I once more became an orphan. I not only  lost my mathematical father, but also the one and only person 
that I can converse with about music, literature,languages, arts and sciences including history, theology,
philosophy and metaphysics, namely  almost everything!\footnote{For a short autobiography one may refer to\\ http://sertoz.bilkent.edu.tr/depo/CemTezerHayatim.pdf. \\
 69.} It should be obvious to the reader that many other share similar feelings with me, \cite{27}.\footnote{For a sincere interview with Cem Tezer: Sertöz, S. Matematiğin Aydınlık Dünyası, TÜBİTAK Popüler Bilim Kitapları, Ankara, 1996. pp. 69.}

\section{His Contributions to mathematics}

Dynamical systems, ergodic theory and global analysis, not surprisingly,  play a center role in his research, \cite{2}--\cite{9} and \cite{23}. 
His thesis {\em "Generalised Solenoids of the Dantzig Type”},  provided the motivation and foundation for some of his later work, \cite{2}.
In particular for his papers in Quarterly  Journal of Mathematics 1989, \cite{4} and in the Mathematische Zeitschrift 1992, \cite{6}. 
Tezer's approach in these papers enabled others to obtain further results in this area.  Under the heading of "Smale spaces" the dynamics on solenoids as considered by Professor Tezer has become by now a prominent  topic. Tezer's findings can  be seen as an early contribution to this line of research.  
He also interested in differential geometry, \cite{10} and \cite{22}, classical mechanics \cite{24} and geometry \cite{15}--\cite{21}.  

In his research articles, he made used of  his own rigorous methods and 
that he developed through reinterpreting the known techniques and supported by a meticulously chosen notation.
In his proofs he reflected the mathematical beauty, neatness and elegance. The articles he published in this sense
was complete, irreproachable and refulgent. It was also his fastidiousness 
and perfectionism that he  followed the masters attitude towards publishing a few but ripe. Being his student 
I knew that he  kept and hid out  many elegantly written unpublished work that were not able to jump his  thresholds.
This fact indeed makes me sad. The only thing  that I can do  at this stage of my life is to  hope and
to dream that one day they will be caught by an appreciative mathematician.

Since his research interest are diverse from dynamical systems and differential geometry to classical mechanics  and
elementary geometry it is almost impossible to cover all of his research, in detail. Instead, I will try to emphasise
the first two area and leave the others, hopefully, as a part of separate research article.

\subsection{Dynamical Systems, ergodic theory and global analysis}

A dynamical system,in its most general form, is a topological space $X$ with a continuous map $f: X \longrightarrow X$ on $X$.
In the category of dynamical systems, characterisation depends on topological conjugacy:

\begin{definition}
\label{def1}
Let $(X,f)$ and $(Y, g)$ be dynamical systems. If there exists a homeomorphism $h: X \longrightarrow Y $
such that $h \circ f = g \circ h$, the given  dynamical systems are called topologically conjugate.
\end{definition}

In theory of dynamical systems "shift equivalence" also takes an important part:

\begin{definition}
\label{def2}
Given the dynamical systems $(X,f)$ and $(Y, g)$  if there exists continuous maps
$f: X \longrightarrow Y$ and $\psi: Y \longrightarrow X $ such that
\begin{equation*}
\phi \circ f = g \circ \phi, ~~  f \circ \psi = \psi \circ g, ~~~ \psi \circ \phi =  f^n, ~~\phi \circ \psi  = g^n
\end{equation*}
\noindent for some $n \in \mathbb{Z}^{+}$ then the two dynamical systems are called shift equivalent.
\end{definition}

Both definition \ref{def1} and definition \ref{def2} may be adapted  to any category.
In the  category of groups, applications of these concepts related to the fundamental groups of 
topological spaces play an important role in reducing the problems associated to topological and differentiable 
maps to the algebraic ones. Another  categoric technique used in topological or differential dynamical systems
is the inverse limit, \cite{3}. The importance of the dynamical system $\displaystyle (\Sigma_f, \sigma_f) $
which is obtained from the topological space 

\[\Sigma_f = \{ (x_i)_{i \in \mathbb Z^{+}} : x_i \in X, f(x_{i + 1}) = x_i,~i \in \mathbb{Z}^{+} \} \]

\noindent and the homemomorphism 
\begin{equation*}
\sigma_f       : \Sigma_f  \longrightarrow  \Sigma_f , ~~  \sigma_f (x_i) = (x_{i+1})
\end{equation*}
\noindent can be seen in the applications that appear in Tezer's respective publications, \cite{3}--\cite{8}, \cite{23}.
In particular, his article entitled as "The shift on the inverse limit of a covering projection", \cite{3} takes the
central position among the others. In the following the problems that appear in different articles of Cem Tezer  
will be investigated under two categories, namely those related to the topological and to the differentiable dynamical systems, respectively.

\subsubsection{His contributions to topological dynamical systems}

Given topological spaces $Z_1, ~ Z_2$ and points $z_i \in Z_i,  ~~ i = 1,2 $, the fundamental group with the base point  
$z_i$ is denoted by $\displaystyle\pi_1  \left(Z_i, Z_i \right), ~~ i = 1,2 $. If $f: Z_1 \longrightarrow Z_2$ is a 
continuous map with $f(z_1) = f(z_2)$, the group homomorphism between $\pi_1  \left(Z_1, Z_1 \right)$ and 
$\pi_1  \left(Z_2, Z_2 \right)$ is defined as 
 \begin{eqnarray*}
  f_{*}                : \pi_1  \left(Z_1, Z_1 \right) &\longrightarrow & \pi_1  \left(Z_2, Z_2 \right)\\
                                           f_{*} ([\gamma])    & =& [f \circ \gamma].
 \end{eqnarray*}
 
Let $X$ and $Y$ be connected, locally path connected, semi-locally simply connected metric spaces,
let $a: X \longrightarrow X $ and $b: Y \longrightarrow Y $  be covering maps and  let $x_0, ~y_0$ be the 
fixed points of the maps $a, ~b$. In \cite{3}, the topological equivalence  between 
$(\Sigma_a, \sigma_a)$ and $(\Sigma_b, \sigma_b)$  is identified with
the equivalence between $\left( \pi_1(X,x_0), a_{*} \right)$ and $\left( \pi_1(Y,y_0), b_{*} \right)$ in the category of groups:

\begin{itemize}
\item[\bf{i.}] If If there exists a homeomorphism $ f: \Sigma_a \longrightarrow \Sigma_b  $ such that 
$f(\xi_o) = \eta_o $  and $f \circ \sigma_a  = \sigma_b \circ f$ then
\[ a_{*} : \pi_1 (X, xo)  \longrightarrow  \pi_1(X, xo) \] 
\noindent  is shift equivalent to
\[ b_{*} : \pi_1 (Y, Yo)  \longrightarrow  \pi_1 (Y, Yo),\] \noindent {\bf Thm 4.1}, \cite{3}.\\
\item[\bf{ii.}] If $\sigma_a$ is topologically conjugate to $\sigma_b$ then
\[ a_{*}: \pi_1 (X, xo)  \longrightarrow  \pi_1(X, xo)  \]
\noindent weakly shift equivalent to 
\[b_{*} : \pi_1 (Y, Yo)  \longrightarrow  \pi_1 (Y, Yo),\] \noindent {\bf Thm 4.2}, \cite{3}.
\end{itemize} 

Weak shift equivalence formalised  in \cite{3} will be defined as follows:

\begin{definition} Group endomorphisms $a : G \longrightarrow G$ and $b: H \longrightarrow H$ are
said to be weakly shift equivalent if there exist $g \in G$ and $h \in H$ such that
$Ad[g] \circ a $  is shift equivalent to $ Ad[h] \circ b$, \cite{8}.
\end{definition}

Cem Tezer also showed that when a shift equivalence between $( \pi_1 (X, xo),a_{*} ) $ and 
$( \pi_1 (Y, Yo) ,b_{*} ) $ is given if there exist a covering map $F: X \longrightarrow Y $ satisfying 
specific conditions then the dynamical  systems  $(\Sigma_a, \sigma_a)$ and  $(\Sigma_b, \sigma_b)$ are topologically equivalent,
{\bf Thm. 4.3}, \cite{3}.

\subsubsection{His contributions to differentiable dynamical systems}

Differentiable dynamical systems that Cem Tezer took into account consists of  a smooth manifold $X$
and an expanding map $a: X \longrightarrow X$.

\begin{definition}
Let $X$ be a differentiable manifold on which there exists a Riemannian $<.~ ,~ . >$ metric. $a: X \longrightarrow X $ will be called
and expanding mapping  if there exists constant $C > 0 $, $\lambda > 1$ 
such that 
\[ || Ta^n(x)(\nu) ||_{a^n(x)} > C \lambda^n || \nu(x)||_x \]
\noindent for any $x \in X$, $\nu \in T_x X $ and $n \in \mathbb Z^{+}$ where $|| \nu(x)||_x $ stands for $<\nu,\nu>_x$, \cite{28}.
\end{definition} 

Let $X$ be a compact manifold and $a: X \longrightarrow X$ be an expanding mapping then the dynamical system
$(\Sigma_a, \sigma_a)$ is called an expanding attractor. It is very well known that if  $a: X \longrightarrow X$ be an expanding mapping then it is a covering map and and  $a$ has fixed points. These properties are frequently made used of  in \cite{3} and \cite{5}.

Cem Tezer, as in the case of covering maps has found 
a necessary and sufficient condition describing the relationship between the topological equivalence  
and the shift equivalence:

\begin{itemize}
\item[\bf{i.}] Let $a: X \longrightarrow X$ and $b: Y \longrightarrow Y$  be expanding maps and let $x_0$
and $y_0$ be fixed points of $\sigma_a$ and $\sigma_b$ respectively. $a_{*}: \pi_1(X,x_0) \longrightarrow  \pi_1(X,x_0)$ 
is shift equivalent to  $b_{*}: \pi_1(Y,y_0) \longrightarrow  \pi_1(Y,y_0) $ iff there exists a homeomorphism 
$F: \Sigma_a \longrightarrow \Sigma_b $,   {\bf Thm. 4.6}, \cite{3}.\\
\item[\bf{ii.}] Let $a: X \longrightarrow X$ and $b: Y \longrightarrow Y$  be expanding maps and let $x_0$
and $y_0$ be fixed points of $\sigma_a$ and $\sigma_b$ respectively. $\sigma_a$
is topologically equivalent to $\sigma_b$ iff $a_{*}: \pi_1(X,x_0) \longrightarrow  \pi_1(X,x_0)$ is weakly shift equivalent to $b_{*}: \pi_1(Y,y_0) \longrightarrow  \pi_1(Y,y_0) $, {\bf Thm. 4.6}, \cite{3}.
\end{itemize}
Another interesting result obtained by Cem Tezer related to the differentiable dynamical systems is the following:
Let $X$ be a compact manifold and $a: X \longrightarrow X$ be an expanding map then the fixed points of this map 
is equal to the Reidemeister number of the homomorphism $a_{*}: \pi_1(X,x_0) \longrightarrow  \pi_1(X,x_0) $. {\bf Prop. 1}, \cite{5}. The algebraic  methods that are introduced in the article are then employed to calculate the number of fixed points and the Artin-Mazur zeta function of the expanding maps of the Klein bottle.

It is an old conjecture that Anosov diffeomorphisms are topologically conjugate to infranilmanifold
automorphisms. Cem Tezer in his  survey article entitled as "Recent geometric developments in the theory 
of Anosov diffeomorfisms" discusses  about the work in this direction which involve the geometric
artifact of connections, \cite{9}. In fact, he gives a definition of connections of class $C^0$ and enumerates some results of  his own concerning the existence of some $C^0$(canonical) connection having some properties {\bf Thm. 4.2, Thm. 4.3} and asserts that it will play a crucial role in the theory of Anosov diffeomorphisms,{Conjecture 2.5}, \cite{9}. Unfortunately,  his life was not
long enough to continue his research in this direction, we may only hope that his intuition will be proved soon.

\subsection{His contributions to differential geometry}

Cem Tezer interested in and worked actively in differential geometry, as well. In particular, he worked not only conformal and invariant structures but also on the geometry of the tangent bundle of a Riemannian manifold.

In his paper "Sur les transformations conformes de la sphère admettant une connexion invariante"  he shows  that a conformal transformation $f$ of the standard sphere $S^n, ~   n = 1, 3, 7$ admits an invariant connection if and only if it admits an invariant
Riemannian metric, {\bf Thm. 2}. As a corollary  he shows that If $n$ is even, the set of conformal transformations of 
$S^n$ which do not admit invariant connections contains an open and dense part of the group
of conformal transformations, \cite{10}.

The paper entitled as "Conformal vector fields with respect to the Sasaki metric tensor" is actually
consists of the results  of the   Ph. D.  thesis with the same title of  his student Fatma Muazzez Şimşir, (the author of this article). The idea of lifting the metric that defines the Riemannian manifold to its tangent bundle to get
a natural metric on its tangent bundle is a very well known idea. The most basic infinitesimal transformations on a Riemannian manifold are Killing vector fields and the conformal vector fields. The tangent bundle of a Riemannian manifold is
itself a manifold  with a natural Riemannian manifold structure induced by the Sasaki
metric tensor field  whereof the Killing vector fields have been completely determined
by S. Tanno.  Let $TM$ be the tangent bundle of a Riemannian manifold $(M, g)$ endowed with the
Sasaki metric $G$ defined by g. In this paper, conformal vector
fields on the tangent bundle $TM$ of ($M, g)$ of dimension at least three with respect to
the Sasaki metric $G$ is characterised, {\bf Thm 1}, \cite{22}. It is  also shown that if $(M, g)$ is a compact Riemannian manifold with dim $ M \geq 3$, then a vector field $A$ is conformal with respect to the Sasaki metric G if and
only if $ A $ is a Killing vector field with respect to G, {\bf Thm 2}, \cite{22}. 

\section{Epilogue and acknowledgements}

When news of his death reached me, I already felt a keen sense of loss but would only begin to understand its full extent during the weeks and months that lay ahead. With academy in the world becoming more aggressive, with academicians exposed to increasing demand of publications without considering the quality , the values that Cem Tezer espoused seem all the more attractive. And as so much of the world is in turmoil from publish and perish fashion of various sorts, including to carrying out research without realising the beauty behind, Cem Tezer's attitude towards was almost unique and from my standpoint the mathematical community  should urgently at
least comprehend and recognise his ideas. 

In general, it is very hard to summarise  someone's  whole life and endeavours in a short article. When it comes to my 
dear supervisor, my mathematical father Cem Tezer even it becomes a tough and daring action for me. I expect that this humble article 
becomes an opportunity  to those who knows him personally to remember him with love, joy and gratitude; and to the others
to acquire some information about his life, his research and his understandings about mathematics.


{\em Bezm-i ezelde tekrar mülaki olabilmek niyazlarımla... Ruhun şad olsun Cem Hocam !}

\end{document}